\magnification=1200
\baselineskip=14pt plus 1pt minus 1pt
\tolerance=1000\hfuzz=1pt 


\font\bigfont=cmr10 scaled\magstep3
\font\secfont=cmr10 scaled\magstep1

\def\section#1#2{\vskip32pt plus4pt \goodbreak \noindent{\bf\secfont#1. #2}
        \xdef\currentsec{#1} \global\eqnum=0 \global\thmnum=0}
\def\appendix#1#2{\vskip32pt plus4pt \goodbreak \noindent
                {\bf\secfont{Appendix #1. }#2}
        \xdef\currentsec{#1} \global\eqnum=0 \global\thmnum=0}

\def\leftheader{}
\def\rightheader{}
\def\header{\headline={\tenrm\ifnum\pageno>1 
                        \ifodd\pageno\rightheader
                         \else\noindent\leftheader 
                          \fi\else\hfil\fi}}

\newcount\thmnum
\global\thmnum=0
\def\prop#1#2{\global\advance\thmnum by 1
        \xdef#1{Proposition \currentsec.\the\thmnum}
        \bigbreak\noindent{\bf Proposition \currentsec.\the\thmnum.}
        {\it#2} }
\def\defin#1#2{\global\advance\thmnum by 1
        \xdef#1{Definition \currentsec.\the\thmnum}
        \bigbreak\noindent{\bf Definition \currentsec.\the\thmnum.} 
        {#2} }
\def\lemma#1#2{\global\advance\thmnum by 1
        \xdef#1{Lemma \currentsec.\the\thmnum}
        \bigbreak\noindent{\bf Lemma \currentsec.\the\thmnum.} 
        {\it#2} }
\def\thm#1#2{\global\advance\thmnum by 1
        \xdef#1{Theorem \currentsec.\the\thmnum}
        \bigbreak\noindent{\bf Theorem \currentsec.\the\thmnum.} 
        {\it#2} }
\def\cor#1#2{\global\advance\thmnum by 1
        \xdef#1{Corollary \currentsec.\the\thmnum}
        \bigbreak\noindent{\bf Corollary \currentsec.\the\thmnum.} 
        {\it#2} }
\def\conj#1#2{\global\advance\thmnum by 1
        \xdef#1{Conjecture \currentsec.\the\thmnum}
        \bigbreak\noindent{\bf Conjecture \currentsec.\the\thmnum.} 
        {\it#2} }

\def\proof{\medskip\noindent{\it Proof. }}

\newcount\eqnum
\global\eqnum=0
\def\num{\global\advance\eqnum by 1
        \eqno({\rm\currentsec}.\the\eqnum)}
\def\eqalignnum{\global\advance\eqnum by 1
        ({\rm\currentsec}.\the\eqnum)}
\def\ref#1{\num  \xdef#1{(\currentsec.\the\eqnum)}}
\def\eqalignref#1{\eqalignnum  \xdef#1{(\currentsec.\the\eqnum)}}

\def\title#1{\centerline{\bf\bigfont#1}}

\newcount\subnum
\def\Alph#1{\ifcase#1\or A\or B\or C\or D\or E\or F\or G\or H\fi}
\def\subsec{\global\advance\subnum by 1 
        \vskip12pt plus4pt \goodbreak \noindent
        {\bf \currentsec.\Alph\subnum.}  }

\def\today{\ifcase\month\or January\or February\or March\or 
        April\or May\or June\or July\or August\or September\or 
        October\or November\or December\fi\space\number\day, 
        \number\year}

\input amssym.def
\input amssym.tex

\def\bR{{\Bbb R}}

\def\bZ{{\Bbb Z}}

\def\cA{{\cal A}}

\def\fA{{\frak A}}

\chardef\o="1C
\def\cK{{\cal K}}
\def\fB{{\frak B}}

\header
\def\leftheader{\centerline{S. KLIMEK}}
\def\rightheader{\centerline{QUANTUM DYNAMICS}}

\def\intsec{I}
\def\ergsyssec{II}
\def\classec{III}
\def\equisec{IV}
\def\reprsec{V}
\def\examplsec{VI}

\def\mapright#1{\smash{\mathop{\longrightarrow}\limits^{#1}}}

\hsize=16.0truecm\vsize=22.5truecm\vglue6.3truecm

{\baselineskip=12pt
\title{Quantum Dynamical Systems}
\vskip 1cm
\title{with Quasi--Discrete Spectrum}
\vskip 1in
\centerline{\bf S\l awomir Klimek}
\footnote{}{Supported in part by the National Science Foundation under
grant DMS--9801612}
\vskip 12pt
\centerline{Department of Mathematics}
\centerline{IUPUI}
\centerline{Indianapolis, IN 46205, USA}
\vskip 1in\noindent
{\bf Abstract.} We study totally ergodic quantum dynamical systems with
quasi--discrete spectrum. We investigate the classification problem
for such systems in terms of algebraic invariants. The results are noncommutative
analogs of (a part of) the theory of Abramov.
\vfill\eject}

\section\intsec{Introduction}
\medskip

Let $(X,\mu)$, $\mu(X)=1$, be a standard Lebesgue space and let
$\alpha:X\rightarrow X$ be an automorphism of $(X,\mu)$.
Then $\alpha$ defines an unitary operator, called the Koopman operator [K],
in $L^2(X,d\mu)$ and denoted by the same letter.

In the important papers [VN] and [HvN], von Neumann and
Halmos classified all classical ergodic systems for which the Koopman
operator has purely discrete spectrum. 
The main result
of their analysis is that such systems are classified by the spectrum, which
forms a discrete subgroup of $U(1)$, and each such a system is conjugate to a
shift on a compact abelian group, the Pontriagin dual of the spectrum. 
Here, and throughout the paper, $U(1)$ is the group of complex numbers with
absolute value 1 and discrete topology.
For a clear account of that result, see e.g. [CFS],  [W], or [Si].

This theory was extended to noncommutative setting by Olsen, Pedersen
and Takesaki [OPT]. It turns out that noncommutative ergodic systems with
discrete spectrum are classified by the spectrum of the automorphism,
which as above is a discrete subgroup $H$ of $U(1)$ and a second cohomology
class of $H$. This theorem is stated more carefully in Section \ergsyssec.

The notions of quasi-eigenvalue and quasi-eigenfunction were introduced by von
Neumann and Halmos [H]. They proved, using those concepts, that there exist
spectrally equivalent but not conjugate automorphisms with mixed spectrum.
Later Abramov [Ab] gave a complete classification of totally ergodic systems
with quasi-discrete spectrum. A topological version of Abramov's
theory for minimal systems was discussed in [HaP], [HoP].

Let us shortly describe what quasi-eigenvalues and quasi-eigenfunctions
are and state the Abramov's theorem. With the above notation $\alpha$ is called
totally ergodic if $\alpha^n$ is ergodic for every $n=1,2,\ldots$. 
Ordinary eigenvectors and eigenvalues of $\alpha$ are called, correspondingly,
quasi-eigenvectors and quasi-eigenvalues of the first order.
A function $f\in L^2(X,d\mu)$
is called a quasi-eigenvector of the second order if
$$
\alpha(f)=\phi f,
$$
where $\phi$ is a quasi-eigenvectors of the first order (i.e.
an eigenvector) of $\alpha$. In such a case $\phi$ is
called a quasi-eigenvalue of the second order. Continuing this process one
obtains quasi-eigenvectors and quasi-eigenvalues of arbitrary order - see
Section \ergsyssec\ for a more precise definition. 
The crucial observation is that, if $\alpha$ is totally ergodic,
quasi-eigenvectors corresponding to different quasi-eigenvalues are
orthogonal. One considers then the situation when
$L^2(X,d\mu)$ has a basis consisting of quasi-eigenvectors of $\alpha$
possibly of arbitrary order.
If this is the case, then we say that $\alpha$ has purely quasi-discrete spectrum.
The Abramov's theorem can be formulated as follows.

\thm\abramthm{{\bf [Ab]}\hskip 5mm
There is a one-to-one correspondence between the conjugacy
classes of totally ergodic dynamical systems with purely quasi-discrete spectrum
and the equivalence classes of pairs $(H,R)$ where H is a discrete abelian group
of the form $H=\bigcup_{n=1}^\infty H_n$ where $H_1\subset H_2\subset\ldots$ is
an increasing sequence of discrete abelian groups, $H_1\subset U(1)$ and $H_1$
has no non-trivial elements of finite order, and $R$ is a homomorphism of $H$
such that for every $n=1,2,\ldots$ the kernel of $R^n$ is the group $H_n$.
}
\medskip

This paper contains an attempt to extend the Abramov theorem to the
quantum mechanical context i.e. when the space $X$ is replaced by a 
noncommutative von Neumann algebra. One case of this program that we were
able to understand fairly completely is when the second order quasi-eigenvectors form a basis in the corresponding $L^2$-space. This assumption is satisfied in the original example that has motivated our work on the subject.
The main results of the paper, Equivalence Theorem and 
Representation Theorem, show that such systems are classified by
quadruples $(H_1,H_2,[r],k)$, called quantum quasi-spectra, where
$H_1$ and $H_2$ are groups, $k:H_2\mapsto H_2$ is an isomorphism and
$[r]$ is (essentially) a $k$ invariant second cohomology class of $H_2$.
Such quadruples are also required to satisfy a number of conditions described 
in Section \equisec. 

It seems that the classification problem in full generality leads to an excessively complicated system of algebraic invariants and is left for
future investigation. In what follows we present a detailed account of the classification theory under the above mentioned additional assumption. 

Our proofs and organization of
the material follow closely that of Abramov's with several important
differences. Among them are:

$\bullet$ The set of quasi-eigenvalues forms a group but not with respect
to operator multiplication but rather a twisted version of it denoted by $*$
in this paper.

$\bullet$ We introduce a natural concept of a normalized basis of
quasi-eigenvectors which simplifies proofs of the Equivalence Theorem and the
Representation Theorem.  

The paper is organized as follows. In Section
\ergsyssec\ we introduce a fairly general setup and precisely formulate 
the problem.  In Section \classec\  we show how to construct group-theoretic invariants for totally
ergodic quantum dynamical systems with purely quasi-discrete spectrum
(of the second order). We prove the equivalence theorem in Section
\equisec, and the representation theorem in Section \reprsec.
Finally, Section \examplsec\ contains a simple example of such a
quantum dynamical system.

\section\ergsyssec{Quantum Ergodic Systems}
\medskip

We begin by reviewing the basic concepts which are used throughout the
paper. We will work within the von Neumann algebra framework, see e.g. [BR],
as this is the natural setup for noncommutative (quantum) ergodic theory. We will
adopt the following definition of a quantum dynamical system.

\defin\qdsdef{
A quantum dynamical system is a quadruple $\left( \fA,G,\alpha ,\tau
\right) $ with the following properties:
\item{(i)}  {$\fA$ is a von Neumann algebra with a separable predual.}
\item{(ii)}  {$G$ is a locally compact abelian group.}
\item{(iii)}  {$\alpha :G\rightarrow \rm{Aut}\left( \fA\right) $ is
an action of $G$ on $\fA$ by von Neumann algebra automorphisms.}
\item{(iv)}  {$\tau $ is a $G$-invariant, normal, faithful state on $\fA$.}
}
\medskip

Since locally compact abelian groups
are amenable, it allows one to define the time average of an observable and
prove ergodic theorems, see e.g. [L], [J], and references
therein. The most relevant are
the groups $G=\bZ$ (in which case the system is called a quantum map)
and $G=\bR$ (in which case the system is called a quantum flow). 

We will denote by $\cK=L^{2}\left( \fA,\tau \right) $ the GNS
representation space of $\fA$ associated with the state $\tau$. 
Since $\fA$ has a separable predual, 
$\cK$ is a separable Hilbert space. It is natural to think of $\cK$ as a
quantum version of the classical Koopman space. The automorphisms $\alpha _{g}$
extend to unitary operators of the $\cK$-spaces. By a slight
abuse of notation, we continue to denote them by $\alpha _{g}$.

\defin\conjdef{
Two quantum dynamical systems $\left( \fA ,G ,
\alpha ,\tau \right) $ and $\left( \fB,G,\beta ,\omega \right) $
are conjugate if there exists an isomorphism of von Neumann algebras 
$\Phi :\fA\rightarrow \fB$ such that
\item{(i)}  {$\Phi \circ \alpha =\beta \circ \Phi $;}
\item{(ii)}  {$\omega \circ \Phi =\tau $.}
}
\medskip

A non-zero element $U\in \cK$ is an eigenvector of $\alpha$ if for every $g\in G$ we have
$\alpha _{g}(U)=\lambda (g)U$, where $\lambda (g)\in U(1)$. Clearly, each 
$g\rightarrow \lambda (g)$ is a character of the group $G$. The set 
Spec$_{p}\left( \alpha \right) $ of all such characters is called the
point spectrum of $\alpha$.

\defin\purdisdef{
A quantum dynamical system $\left( \fA,G,\alpha ,\tau \right) $ is
called a system with purely discrete spectrum if $\cK$ has an orthonormal basis
consisting of eigenvectors of $\alpha $. }
\medskip

As a consequence of the separability
assumption, Spec$_p\left(\alpha\right)$ is a countable subset of the dual group
$\widehat{G}$.

Ergodic theory of von Neumann algebras has been studied by many authors. For
references and a variety of results, see e.g. [C], [KL1,2], [KLMR], [L] and [J].
For our purposes, the following definition of quantum ergodicity will be
sufficient.

\defin\ergdef{
A quantum dynamical system $\left( \fA,G,\alpha ,\tau \right) $ is
called ergodic if the only $G$-invariant elements of $\cK$
are scalar multiples of $I$.
}
\medskip

Equivalently, the joint eigenspace of $\alpha _{g}$'s corresponding to the eigenvalue $1$ is one dimensional and consists of the scalar multiples of
the identity operator. For quantum ergodic systems, the time and ensemble
averages of an observable are equal. Also one has the following classification
theorem due to Olsen, Pedersen and Takesaki [OPT].

\thm\optthm{{\bf [OPT]}\hskip 5mm
There is a one-to-one correspondence between the conjugacy
classes of ergodic quantum dynamical systems with purely discrete spectrum
and the family of pairs $(H,\sigma)$ where $H\subset\widehat{G}$ is a
discrete group and $\sigma$ is a second cohomology class of $H$.
}
\medskip

In fact, in analogy with the commutative theory, every quantum dynamical system
is conjugate to a shift on the noncommutative deformation of $\widehat{H}$
determined by $\sigma$ - see {\bf [OPT]}.

\defin\totergdef{
A quantum dynamical system $\left( \fA,G,\alpha ,\tau \right) $ is
called totally ergodic if for every $g\in G$ individually, the only elements of 
$\cK$ invariant under $\alpha_g$, are scalar multiples
of $I$. 
}
\medskip

For an example of ergodic but not totally ergodic quantum dynamical system
see Section \examplsec.

We shall call the eigenvectors of $\alpha$ quasi-eigenvectors of the first order.
Similarly, eigenvalues of $\alpha$ are called quasi-eigenvalues of the first
order. The set of normalized quasi-eigenvectors of the first order is denoted by
$G_1$ while the set of all quasi-eigenvalues of the first order is denoted by
$H_1$. We define the set $G_n$ of normalized quasi-eigenvectors of $n$-th order
and the set $H_n$ of quasi-eigenvalues of $n$-th order inductively. Suppose that
$G_n$ and $H_n$ are defined.

\defin\quasieigdef{With the above notation, a non-zero element $U\in \cK$ is called a quasi-eigenvector of order $n+1$ of $\alpha$ if 
$\alpha _{g}(U)=\lambda (g)U$, where $\lambda (g)\in \fA\cap G_n$.
Then $\lambda$ is called a quasi-eigenvalue of order $n+1$.
}
\medskip

\defin\purdisdef{
A quantum dynamical system $\left( \fA,G,\alpha ,\tau \right) $ is
called a system with purely quasi-discrete spectrum if $\cK$ has an orthonormal
basis consisting of quasi-eigenvectors of $\alpha $ of possibly arbitrary
orders.}
\medskip

The subject of this paper is the classification problem for (noncommutative)
totally ergodic systems with quasi-discrete spectrum. This is to be solved by
constructing a complete set of algebraic invariants of such systems.

\section\classec{Classification of Quasi-Discrete Systems}
\medskip

\noindent In this paper we tackle the program described in the previous section
under the following additional assumptions:

\item{1.} We consider only $G=\bZ$, i.e. quantum maps. The automorphism
$\alpha_1$ corresponding to the generator 1 of $\bZ$ will simply be denoted by
$\alpha$.
\item{2.} We asume that $\cK$ has an orthonormal basis consisting of
the second order quasi-eigenvectors of $\alpha $.
\item{3.} We require that $\tau$ is a normalized trace.

\noindent Additionally, throughout the rest of the paper we assume that the
system $\left( \fA,\bZ,\alpha,\tau\right)$ is ergodic. We do explicitly
mention when total ergodicity is used.

With extra effort the classification program can be presumably carried out
for arbitrary abelian locally compact groups and, what is most challenging,
arbitrary quasi-discrete spectrum. The trace assumption is used in the proof of unitarity in the following proposition and possibly is not really needed. In any case it seems likely that ergodicity and discreteness of the quasi-spectrum will force any invariant state to be a trace.

Every constant is an eigenvector belonging to the eigenvalue
$\lambda=1$, and therefore $H_1\subset G_1$. Moreover, obviously:
$$
H_1\subset H_2\subset G_1\subset G_2. \ref{\inclref}
$$

\prop\unitaryprop{Let $\lambda$ be an eigenvalue of $\alpha$. If $U_\lambda\in\cK$
is a normalized second order quasi-eigenvector of $\alpha$: 
$$
\alpha\left(U_\lambda\right)=\lambda U_\lambda\ ,\ref{\alpharef}
$$
then $U_\lambda\in\fA$ and $U_\lambda$ is unitary.
}
\proof This needs a little von Neumann algebras theory from [Ar].
Let $P^\natural\subset L^{2}\left(\fA,\tau\right)$ be the closure of
$\fA_+ 1$, where  $\fA_+$
is the positive part of $\fA$ and where 
$1\in\fA\subset L^{2}\left(\fA,\tau\right)$ 
is the unit in $\fA$.
It follows from this definition that $P^\natural$ is invariant under $\alpha$.
It is known that every $x\in L^{2}\left(\fA,\tau\right)$ has a unique
decomposition:
$$
x=u\,|x|,
$$
where $u\in\fA$ is a partial isometry and $|x|\in P^\natural$. 
Write $U_\lambda=u\,|U_\lambda|$ in \alpharef . Then:
$$
\alpha(u)\alpha(|U_\lambda|)=\left(\lambda u\right)|U_\lambda|
$$
It follows that
$|U_\lambda|$ is an invariant vector for $\alpha$ and so, by ergodicity,
it is equal to 1. But that means that $U_\lambda\in\fA$. Applying the
ergodicity assumption to $U_\lambda^* U_\lambda$ we see that 
$U_\lambda^*U_\lambda=1$.

Since $1-U_\lambda U_\lambda^*$ is positive and 
$\tau\left(1-U_\lambda U_\lambda^*\right)=
\tau\left(1-U_\lambda^* U_\lambda\right)=0$
we see that $U_\lambda$ is unitary.
$\square$

\prop\simplicityprop{If $U,V\in G_2$ belong to the same quasi-eigenvalue 
$\lambda$ then there is a constant $C$, $|C|=1$, such that $U=CV$.
}
\proof Applying $\alpha$ to $U^{-1}V$ yields:
$$
\alpha(U^{-1}V)=U^{-1}\lambda^{-1}\lambda V=U^{-1}V.
$$
It follows from ergodicity of $\alpha$ that $U^{-1}V$ is a constant.
$\square$

Let us recall from [OPT] the following structural result about $G_1$.

\prop\ordeigenprop{For each pair $\lambda ,\mu \in H_1$, we have
$$
U_{\lambda }U_{\mu }=\sigma \left( \lambda ,\mu \right) U_{\mu }U_{\lambda },
\ref{\comrel}
$$
where $U_{\lambda },U_{\mu }\in G_1$ are the corresponding eigenvectors and
$\sigma :H_1 \times H_1 \rightarrow U\left( 1\right) $.  Furthermore,
$\sigma$ has the following properties: 
$$
\sigma \left( \lambda ,\lambda \right) =1,  \ref{\normaliz}
$$
$$
\sigma \left( \lambda ,\mu \nu \right) =\sigma \left( \lambda ,\mu \right)
\sigma \left( \lambda ,\nu \right) ,  \ref{\charprop}
$$
and 
$$
\sigma \left( \mu ,\lambda \right) =\sigma \left( \lambda ,\mu \right) ^{-1}.
\ref{\symmetry}
$$
}
A map $\sigma :H_1 \times H_1 \rightarrow U\left( 1\right) $
satisfying \normaliz, \charprop, \symmetry\ is called a symplectic bicharacter.

The following lemma deals with effects of noncommutativity of $\fA$
on the classification problem.
\lemma\techlemma{
\item{(i)} If $U_\lambda\in G_2$ belongs to quasi-eigenvalue $\lambda\in H_2$
then there exist a number $\phi(\lambda)\in U(1)$ such that
$$
U_\lambda^{-1}\lambda U_\lambda=\phi(\lambda)\lambda.
$$
\item{(ii)} If $U\in G_2$ and $V\in G_1$ then $UVU^{-1}\in G_1$.
}
\proof We verify by direct calculation that $U_\lambda^{-1}\lambda U_\lambda$
and $\lambda$ belong to the same eigenvalue of $\alpha$. Consequently,
\simplicityprop\ implies item $(i)$.

If $U\in G_2$ belongs to $\lambda\in H_2$, $\lambda\in H_2\subset G_1$
belongs to $R(\lambda)\in H_1$, 
and $V\in G_1$ belongs to $\mu\in H_1$, then we compute:
$$\eqalign{
\alpha(UVU^{-1})&=\lambda U\mu V U^{-1}\lambda^{-1}=
{\mu\over \phi(\lambda)}\lambda UV\lambda^{-1}U^{-1}\cr
&=\mu U\lambda V\lambda^{-1}U^{-1}=\mu\sigma(R(\lambda),\mu)
UVU^{-1}\cr
}\ref{\conjformref}
$$
which proves $(ii)$. In the above calculation we used $(i)$ twice
as well as \ordeigenprop .
$\square$

If $\lambda,\mu\in H_2$ and $U_\lambda\in G_2$ is a quasi-eigenvector belonging
to $\lambda$ we define the following product on $H_2$:
$$
\lambda * \mu :=\lambda  U_\lambda \mu U_\lambda^{-1}.\ref{\stardefref}
$$

\prop\groupprop{ Each of the sets $H_1, G_1, G_2$ is a group under operator
multiplication while $H_2$ is a group under $*$ multiplication. Moreover
$H_1\subset H_2$ is a subgroup.
}
\proof The fact that $H_1$ and $G_1$ are groups follows from [OPT]
so we need to concentrate on $H_2$ and $G_2$. 
We first verify that the right hand side of \stardefref\ is in $G_1$:
$$
\alpha(\lambda  U_\lambda \mu U_\lambda^{-1})=
R(\lambda)R(\mu)\sigma(R(\lambda),R(\mu))\cdot 
\lambda  U_\lambda \mu U_\lambda^{-1}\ref{\rofstarref}
$$
by \conjformref .
Here $R(\lambda)$ and $R(\mu)$ are eigenvalues corresponding
to eigenvectors $\lambda$ and $\mu$. 
Additionally:
$$
\alpha(U_\lambda U_\mu)=\lambda U_\lambda\mu U_\mu=
\lambda  U_\lambda \mu U_\lambda^{-1}\cdot U_\lambda U_\mu
=\lambda * \mu\cdot U_\lambda U_\mu\ref{\prodruleref}
$$
so that $\lambda * \mu\in H_2$.
Consequently the $*$- product
is well defined. The identity operator $1\in\fA$ is the unit for
this multiplication.
Since 
$$\alpha(U_\lambda^{-1})={\lambda^{-1}\over \phi(\lambda)}\cdot
U_\lambda^{-1}$$
the $*$ inverse of $\lambda$ is 
$$I(\lambda):={\lambda^{-1}\over \phi(\lambda)}$$
with $\lambda^{-1}$ the operator multiplication inverse. Associativity
of the $*$ multiplication follows from \rofstarref\ 
which also shows that $G_2$ is a group under operator multiplication.
Finally if $\lambda,\mu\in H_1$ then $\lambda*\mu=\lambda\mu$. 
$\square$

We define a map $R: G_2\rightarrow H_2$ by $R(U):=\lambda$ if
$\alpha(U)=\lambda U$. In other words, $R$ assigns to a quasi-eigenvector the
corresponding quasi-eigenvalue. Clearly $R$ maps $G_1\subset G_2$ into
$H_1\subset H_2$. Also $R$ maps $H_2\subset G_1$ into $H_1$.

\prop\rmapprop{The mapping $R:H_2\rightarrow H_1$ has the following properties:
\item{(i)} For every $\lambda\in H_2$ and $\mu\in H_1$ we have
$\mu\,\sigma\left(\mu, R(\lambda)\right)\in H_1$ and
$$
\lambda*\mu*I(\lambda)=\mu\,\sigma\left(R(\lambda),\mu\right).
\ref{\rpropref}
$$
In particular, $H_1$ is a normal subgroup of $H_2$.
\item{(ii)} $R$ is a ``twisted" homomorphism:
$$
R(\lambda*\mu)=R(\lambda)*\lambda*R(\mu)*I(\lambda)
=R(\lambda)R(\mu)\sigma(R(\lambda),R(\mu)).\ref{\twist}
$$
\item{(iii)} The kernel of $R$ is the group $H_1$.
}
\proof Item $(i)$ is just a rephrasing of \conjformref\ 
and item $(ii)$ follows directly from \rofstarref .
Item $(iii)$ is a consequence of ergodicity of $\alpha$,
as eigenvectors corresponding to eigenvalue $\lambda=1$ are proportional
to the identity.
$\square$

Let $N:=$Image of $R\subset H_1$. Equip $N$ with the following
product:
$$
n_1*n_2:=n_1n_2\sigma(n_1,n_2)\in N,
$$
where the last inclusion follows from \rmapprop, item (i). It is easy
to see that $N$ is a group with respect to this product and
$R:H_2\mapsto R$ is a homomorphism. Consequently, we have the following
short exact sequence of groups:
$$
1\quad\mapright{}\quad H_1\quad\mapright{}\quad
H_2\quad\mapright{R}\quad N\quad\mapright{}\quad
1\, .\ref{\exactseq}
$$
This sequence is an extension with abelian kernel, and the $N$-module
structure on $H_1$ is given by \rpropref, see [B].
\prop\countprop{The group $H_2$ is at most countable, and, assuming that
$\alpha$ is totally ergodic, $H_1$ has no nontrivial elements of finite order. }
\proof Since $\alpha$ is assumed to be totally ergodic no nontrivial
elements of finite order in $H_1$ can exist. Also $H_1$ is at most countable
as a consequence of separability of $\cK$. Since $R$ defines a one-to-one map
$H_2/H_1\mapsto H_1$, the group $H_2$ is at most countable.
$\square$

If $U$ belongs to $\lambda\in H_2$ then $\alpha(U)$ belongs to
$R(\lambda)*\lambda$. Thus it makes sense to study the properties of the
map:
$$
k(\lambda):=R(\lambda)*\lambda.\ref{\kdefref}
$$
\prop\kprop{The map $k$ defined by \kdefref\ is an isomorphism
of $H_2$. Moreover $k(\lambda)*I(\lambda)\in H_1$ and
$k(\lambda)=\lambda$ iff $\lambda\in H_1$.
}
\proof $k$ is a homomorphism since
$$\eqalign{
k(\lambda*\mu)&=R(\lambda*\mu)*\lambda*\mu=
R(\lambda)*\lambda*R(\mu)*I(\lambda)*\lambda*\mu\cr
&=R(\lambda)*\lambda*R(\mu)*\mu=k(\lambda)*k(\mu)
}$$
by \rmapprop. The inverse of $k$ is 
$k^{-1}(\lambda)=R(\lambda)^{-1}*\lambda$. Next 
$k(\lambda)*I(\lambda)=R(\lambda)$ so it is in $H_1$. Finally
$k(\lambda)=\lambda$ iff $R(\lambda)=1$ so $\lambda\in H_1$.
$\square$

\prop\orthprop{If the automorphism $\alpha$ is totally ergodic, then
quasi-eigenvectors belonging to different quasi-eigenvalues are orthogonal
in $\cK$. 
}
\proof The statement is true for ordinary eigenvectors. Let $\cK_1$
be the closed subspace of $\cK$ spanned by $G_1$, and let $\cK_2$ be its
orthogonal complement. The assumption of total ergodicity of $\alpha$
is used in the following lemma which says that quasi-eigenvector
which is not an eigenvector can not be a linear combination
of eigenvectors.
\lemma\toterglemma{Suppose $U\in G_2$ is not in $G_1$ and
belongs to $\lambda\in H_2$. Then $U\not\in\cK_1$.
}
\proof Assume that 
$$
U=\sum_{\mu\in H_1}a_\mu U_\mu .\ref{\expanref}
$$
We can compute
$\alpha^n(U)$ in two different ways. First use \expanref\ and apply
$\alpha^n$ to each $U_\mu$. This yields:
$$
\alpha^n(U)=\sum_{\mu\in H_1}a_\mu' U_\mu,
$$
where $a_\mu'$ differs from $a_\mu$ by a phase. Secondly, use
$\alpha(U)=\lambda U$ n-times and then expand:
$$
\alpha^n(U)=\sum_{\mu\in H_1}a_\mu'' U_{R(\lambda)^n\mu},
$$
where, as before, $a_\mu''$ differs from $a_\mu$ by a phase. By \countprop\ 
$R(\lambda)^n\mu$ are all different. Consequently, for any $\mu$ there is
an infinite number of coefficients in \expanref\ equal, up to a phase,
to $a_\mu$, and so they must be zero. 
$\square$

Returning to the proof of \orthprop, if $U\in G_2$ and not in $G_1$, then
we claim that $U$ is in $\cK_2$. In fact, let
$U=U_1+U_2$ be the orthogonal decomposition of $U$ with respect to
$\cK=\cK_1\oplus\cK_2$. It follows from \toterglemma\ that $U_2\not=0$.
Since $\alpha$ is unitary, $\alpha(U_1)\in\cK_1$ and $\alpha(U_2)\in\cK_2$.
Moreover $\lambda U_1\in\cK_1$ because $G_1$ forms a group. For the same reason
$\lambda U_2\in\cK_2$ as:
$$
(\mu,\lambda U_2)=(\lambda^{-1}\mu,U_2)=0,
$$
for $\mu\in G_1$. Consequently we have $\alpha(U_1)=\lambda U_1$
and $\alpha(U_2)=\lambda U_2$ which implies, in view of \simplicityprop,
that $U_1=CU_2$. This can happen only if $C=0$ as $U_1$ and $U_2$ belong
to perpendicular subspaces of $\cK$.

It remains to prove that if $U,V\in G_2$
are not in $G_1$ and belong to different quasi-eigenvalues 
$\lambda,\mu\in H_2$ then $U,V$ are orthogonal. But this is the same as proving
that $U^{-1}V$ is orthogonal to $1\in\cK_1$. Since $G_2$ is a group
with respect to operator multiplication, $U^{-1}V\in G_2$ and belongs
to quasi-eigenvalue $I(\lambda)*\mu$. If $U^{-1}V$ is not in $G_1$
then the orthogonality follows from the previous argument. It remains to
consider the case when $U^{-1}V\in G_1$. But two elements of $G_1$
are orthogonal unless they belong to the same eigenvalue, and, since
$\lambda\not=\mu$, $I(\lambda)*\mu\not= 1$.
$\square$

\cor\statecor{For every $\lambda\in H_2$ we have:
$$
\tau(U_\lambda)=\cases{1&if $\lambda=1$\cr
0&otherwise.\cr
}
$$
}
\proof This is a direct consequence of \orthprop\ and $\tau(U_\lambda)=
(1,U_\lambda)$.
$\square$

\section\equisec{Equivalence Theorem}
\medskip

In this section we spell out the complete set of group theoretic invariants
for totally ergodic quantum dynamical systems with quasi-discrete spectrum of
the second order. The equivalence theorem proved here says that if two
such systems have the same set of invariants then they are conjugate.

If $H$ is a group, then a function $r:H\times H\rightarrow U\left( 1\right) $ is
called a 2-cocycle if 
$$
r\left( \lambda ,\mu \right) r\left( \lambda \mu ,\nu \right) =r\left(
\lambda ,\mu \nu \right) r\left( \mu ,\nu \right) ,  \ref{\cocycleref}
$$
for all $\lambda ,\mu ,\nu \in H$.
A $2$-cocycle 
$r$ is called {\it trivial} if there is a function $d:H\rightarrow U\left(
1\right) $, such that $r\left( \lambda ,\mu \right) =d\left( \lambda \mu
\right) /d\left( \lambda \right) d\left( \mu \right) $.
The set of equivalence classes of 2-cocycles mod trivial 2-cocycles is the second
cohomology group $H^2(H)$ of group $H$ (with values in $U(1)$).

\lemma\ulem{
Let $\left( \frak{A},\bZ,\alpha ,\tau \right) $ be a totally
ergodic quantum dynamical system with purely quasi-discrete spectrum of the
second order. Choose an orthonormal basis $\{U_\lambda\}$, $\lambda\in H_2$,
in $\cK$, consisting of quasi-eigenvalues of $\alpha$ and such that $U_1=1$.
Then for each pair $\lambda,\mu \in H_2$, 
$$
U_{\lambda }U_{\mu }=r\left( \lambda ,\mu \right) U_{\lambda *\mu },
\ref{\multiplicref}
$$
where $r\left( \lambda ,\mu\right)$ is a 2-cocycle on $H_2$. Moreover,
any other orthonormal basis of $\cK$ consisting of quasi-eigenvectors
of $\alpha$ leads to a cohomologous $r$ and $\fA$ is linearly spanned
by $\{U_\lambda\}$.
}
\proof \multiplicref\ is a consequence of \simplicityprop , \prodruleref.
The associativity of the operator multiplication implies that $r$ is
a cocycle. If $\{V_\lambda\}$ is any other orthonormal basis of $\cK$
consisting  of quasi-eigenvectors of $\alpha$ then $V_\lambda=d(\lambda)
U_\lambda$, $d(\lambda)\in U(1)$, and $d(\lambda)$ gives the equivalence of the 
corresponding cocycles. Finally, since $U_\lambda$ is a basis in
$\cK$ it follows that $\cA$ is a $\sigma$-weakly closure of the linear span 
of $\{U_\lambda\}$.
$\square$

Since $H_2\subset G_1$, given a choice of a basis in $\cK$ we can
write for any $\lambda\in H_2$:
$$
\lambda=C(\lambda)U_{R(\lambda)},\ref{\decompref}
$$
where $C(\lambda)\in U(1)$. The main properties of the coefficients
$C(\lambda)$ are summarized in the following lemma.
\lemma\clem{ With the above notation we have:
$$
C(\lambda*\mu)=C(\lambda)C(\mu){r(\lambda,R(\mu))\,
r(R(\lambda),\lambda*R(\mu)*I(\lambda))
\over r(\lambda*R(\mu)*I(\lambda),\lambda)}.\ref{\cpropref}
$$
Additionally, if $\lambda\in H_1$ then $C(\lambda)=\lambda$.
}
\proof Proof is a straightforward calculation using \multiplicref,
\decompref, and \rmapprop\ which we omit.
$\square$

Let $D(\lambda)$ be the following $U(1)$-valued function on $H_2$:
$$
D(\lambda)=\cases{\lambda&if $\lambda\in H_1$\cr
1&otherwise.\cr
}\ref{\ddefref}
$$
We shall show below that one can choose a basis $\{U_\lambda\}$, $\lambda\in H_2$,
in $\cK$, consisting of quasi-eigenvalues of $\alpha$, such that the
matrix elements of $\alpha$ are particularly simple.
\prop\basisprop{ There is a basis $\{U_\lambda\}$, $\lambda\in H_2$,
in $\cK$, consisting of quasi-eigenvalues of $\alpha$, such that
$$
\alpha(U_\lambda)=D(\lambda)U_{k(\lambda)}.\ref{\basisnorm}
$$
Such a basis will be called a {\sl normalized} basis.
}
\proof Notice that \basisnorm\ says that
$\alpha(U_\lambda)=\lambda U_\lambda$ is $\lambda\in H_1$, which is
always true, and $\alpha(U_\lambda)=U_{k(\lambda)}$ if $\lambda\notin H_1$.
Consider the orbits of $k$. If $\lambda\in H_1$ then $k(\lambda)=\lambda$
and $H_1$ is the set of fixed points for $k$. If $\lambda\notin H_1$
then $k^n(\lambda)=R(\lambda)^n*\lambda$ and, as $H_1$ has no elements of
finite order, all $k^n(\lambda)$ are different for different $n\in\bZ$.
Choose one element $s(\lambda)$ from each orbit $k^n(\lambda)$,
so that each $\lambda$ can be uniquely written as $\lambda=k^n(s(\lambda))$.
Choose $U_{s(\lambda)}$ arbitrarily and set
$$
U_\lambda:=\alpha^n\left(U_{s(\lambda)}\right).
$$
Since $U_{k(\lambda)}=\alpha^{n+1}\left(U_{s(\lambda)}\right)$,
\basisnorm\ is clearly satisfied.
$\square$

Let $\{U_\lambda\}$ be a normalized basis and let $r(\lambda,\mu)$ be the
corresponding 2-cocycle on $H_2$. Applying $\alpha$ to \multiplicref\
we infer that
$$
{r(k(\lambda),k(\mu))\over r(\lambda,\mu)}=
{D(\lambda*\mu)\over D(\lambda)D(\mu)}.\ref{\rnormref}
$$
Such a cocycle will be called a {\sl normalized}
cocycle. If $V_\lambda=d(\lambda)U_\lambda$, $d(\lambda)\in U(1)$ is
another normalized basis then
$$
d(k(\lambda))=d(\lambda).\ref{\dkref}
$$
By $H^2_k(H_2)$ we denote the set of equivalence classes of
normalized 2-cocycles on $H_2$ modulo $k$-invariant coboundaries \dkref.

\bigbreak\noindent
{\bf Remark.} If $H_2$ is abelian the set $H^2_k(H_2)$ can be alternatively
described as follows. Let $\tilde D$ be a homomorphism of $H_2$ into $U(1)$
extending the natural embedding $H_1\subset U(1)$. Such an extension is always
possible for abelian groups [Ab]. Then, just like in Proposition \ddefref,
a basis $\tilde U_\lambda$ can be constructed satisfying
$\alpha(\tilde U_\lambda)=\tilde D(\lambda)\tilde U_{k(\lambda)}$. The
corresponding 2-cocycle $\tilde r$ on $H_2$ is then k-invariant by an
analog of \rnormref, and cohomologous to $r$ by \ulem. So, in this case,
$H_k^2(H_2)$ is the second group of k-invariant cohomologies of $H_2$.
In general, when $H_2$ is not necessarily abelian, it is desirable to have a
better description of $H_k^2(H_2)$.
\medskip
Let us denote by $[r]$ the cohomology class of $r$ in $H_k^2(H_2)$.
When restricted to $H_1$ the conditions \rnormref\ and \dkref\ are void.
Moreover, since $H_1$ is abelian, there is a one-to-one
correspondence between the second cohomology classes $[r]$ and
symplectic bicharacters $\sigma$, see \ordeigenprop . The correspondence
is given by:
$$
r\left( \lambda ,\mu \right) =\sigma \left( \lambda ,\mu \right) r\left( \mu
,\lambda \right) ,  \ref{\rsigmarel}
$$
see [OPT].

So far to a totally ergodic system with purely quasi-discrete spectrum of the 
second order we have associated the following algebraic structure:
\item{1.} A countable abelian group $H_1\subset U(1)$ which has no nontrivial
elements of finite order.
\item{2.} A countable group $H_2$, such that $H_1\subset H_2$ is a 
normal subgroup.
\item{3.} An isomorphism $k:H_2\mapsto H_2$ such that 
$k(\lambda)*\lambda^{-1}\in H_1$ and $k(\lambda)=\lambda$ iff $\lambda\in H_1$.
\item{4.} A cohomology class $[r]$ in $H_k^2(H_2)$.

\defin\qssysdef{A quadruple $(H_1,H_2,[r],k)$ satisfying conditions 1-4
above is called a quantum quasi-spectrum.
}
\medskip

\defin\sqsymiso{Two quantum quasi-spectra $(H_1,H_2,[r],k)$ and
$(H_1',H_2',[r'],k')$ are called isomorphic if
\item{(i)} $H_1=H_1'$.
\item{(ii)} There exists an isomorphism $\phi$ of the groups
$H_2$ and $H_2'$ leaving fixed all the elements of the group $H_1=H_1'$ and
such that
$$
k=\phi^{-1}k'\phi,\hskip 5mm [r]=\phi^*[r'],
$$
where $\phi^*$ is the induced isomorphism of the cohomology groups.
}
\medskip

We are now prepared to prove the following theorem which is the main result
of the section.

\thm\equivthm{{\bf (Equivalence Theorem)}\hskip 5mm
Let $\left( \frak{A},\bZ,\alpha ,\tau \right) $ and 
$\left( \frak{B},\bZ,\beta,\omega \right) $ be two totally
ergodic quantum dynamical systems with purely quasi-discrete spectrum of the
second order, and let $(H_1(\alpha),H_2(\alpha),[r_\alpha],k_\alpha)$ and 
$(H_1(\beta),H_2(\beta),[r_\beta],k_\beta)$ denote
the corresponding  quantum quasi-spectra.
The following statements are equivalent:
\item{(i)}  The quantum quasi-spectra 
$(H_1(\alpha),H_2(\alpha),[r_\alpha],k_\alpha)$ and
$(H_1(\beta),H_2(\beta),[r_\beta],k_\beta)$ are isomorphic;
\item{(ii)}  $\left( \frak{A},\bZ,\alpha ,\tau \right) $ and 
$\left( \frak{B},\bZ,\beta,\omega \right) $ are conjugate.
}
\proof Only $(i)\to(ii)$ is non trivial.
Let $\cK(\alpha)$ and $\cK(\beta)$ be the corresponding GNS Hilbert spaces.
We are going to construct a conjugation $\Phi:\frak{A}\mapsto\frak{B}$
as an isomorphism implemented by a unitary map $Q:\cK(\alpha)\mapsto
\cK(\beta)$. Let $\{U_\lambda\}$ and $\{V_\mu\}$ be normalized orthonormal basis
in $\cK(\alpha)$ and $\cK(\beta)$ correspondingly, consisting of
quasi-eigenvectors. Set:
$$
Q\left(U_\lambda\right):=V_{\phi(\lambda)},\ref{\qdefref}
$$
where $\phi$ is an isomorphism of $H_2(\alpha)$ and $H_2(\beta)$.
By \ulem\ we have $U_{\lambda_1}U_{\lambda_2}=r_\alpha(\lambda_1,\lambda_2)
U_{\lambda_1*\lambda_2}$
and $V_{\mu_1}V_{\mu_2}=r_\beta(\mu_1,\mu_2)V_{\mu_1*\mu_2}$. Since
$r_\alpha$ and $\phi^*r_\beta$ are cohomologous, we may assume,
renormalizing $V_\mu$ if necessary, that
$$
r_\alpha(\lambda_1,\lambda_2)=r_\beta\left(\phi(\lambda_1),\phi(\lambda_2)
\right).\ref{\rpropref}
$$
We can deduce from \rpropref\ 
that $\Phi(U_\lambda):=QU_\lambda Q^{-1}=V_{\phi(\lambda)}$ as follows:
$$\eqalign{
QU_{\lambda_1} Q^{-1}V_{\phi(\lambda_2)}&=QU_{\lambda_1}U_{\lambda_2}=
r_\alpha(\lambda_1,\lambda_2)QU_{\lambda_1*\lambda_2}=
r_\alpha(\lambda_1,\lambda_2)V_{\phi(\lambda_1*\lambda_2)}\cr
&=r_\beta\left(\phi(\lambda_1),\phi(\lambda_2)\right)
V_{\phi(\lambda_1)*\phi(\lambda_2)}=
V_{\phi(\lambda_1)}V_{\phi(\lambda_2)}\cr
}$$
But $\fA$ and $\fB$ are linearly generated by, correspondingly, $U_\lambda$
and $V_\mu$ and so $\Phi$ extends to an isomorphism of $\fA$ and $\fB$.
A straightforward calculation verifies that $\Phi \circ \alpha=\beta \circ
\Phi$:
$$\eqalign{
(\Phi \circ \alpha)\,U_\lambda&=D(\lambda)\,\Phi\left(U_{k_\alpha(\lambda)}\right)
=D(\lambda)\,V_{\phi(k_\alpha(\lambda))}
=D(\lambda)\,V_{k_\beta(\phi(\lambda))}\cr
&=\beta\left(V_{\phi(\lambda)}\right)
=(\beta \circ\Phi)\,U_\lambda
}$$
Also $\omega(\Phi(U_\lambda)) =\tau(U_\lambda)$
by \statecor. It follows that $\left( \frak{A},\bZ,\alpha ,\tau \right) $ and 
$\left( \frak{B},\bZ,\beta,\omega \right) $ are conjugate.
$\square$

\section\reprsec{Representation Theorem}
\medskip
In this section we prove a representation theorem which says that
for any system of invariants (i.e. a quantum quasi-spectrum) there is
a corresponding quantum dynamical system with exactly that system
of invariants. Consequently, the correspondence between the conjugacy
classes of totally ergodic systems with purely quasi-discrete spectrum
and the isomorphism classes of quantum quasi-spectra is onto.

\thm\reprthm{{\bf (Representation Theorem)}\hskip 5mm
Let $(H_1,H_2,[r],k)$ be a quantum quasi-spectrum. There exists a totally
ergodic quantum dynamical system $\left( \fA,\bZ, \alpha ,\tau \right) $ 
with purely quasi-discrete spectrum such that
its quantum quasi-spectrum is isomorphic to $(H_1,H_2,[r],k)$
}
\proof 
Consider $\cK:=l^2(H_2)$ and let $\{\phi_\lambda\}$ be the canonical
basis in $\cK$. Define $\fA$ to be the von Neumann algebra generated 
by the following operators $U_\lambda$:
$$
U_\lambda\phi_\mu:=r(\lambda,\mu)\phi_{\lambda*\mu},\ref{\udefref}
$$
where $r(\lambda,\mu)$ is a normalized 2-cocycle on $H_2$ corresponding to $[r]$.
For any $f\in\cK$ we obtain
$$
U_\lambda f(\mu)=r(\lambda,I(\lambda)*\mu)f(I(\lambda)*\mu).
$$
It follows that
$$
U_{\lambda }U_{\mu }=r\left( \lambda ,\mu \right) U_{\lambda *\mu },
$$
Then set
$$
\beta \phi_\lambda:=D(\lambda) \phi_{k(\lambda)},\ref{\betadefref}
$$
where $D(\lambda)\in U(1)$ was defined in \ddefref.
Equivalently, or any $f\in\cK$ we have
$$
\beta f(\lambda)=D\left(k^{-1}(\lambda)\right)
f\left(k^{-1}(\lambda)\right)=D\left(\lambda\right)
f\left(k^{-1}(\lambda)\right),\ref{\betaform}
$$
since $D(\lambda)$ is $k$ invariant.
$\beta$ is a unitary operator in $\cK$ with the inverse given by
$$
\beta^{-1}\phi_\lambda={1\over D\left(k^{-1}(\lambda)\right)}\,
\phi_{k^{-1}(\lambda)},
$$
or, equivalently, for any $f\in\cK$
$$
\beta^{-1}f(\lambda)={1\over D(\lambda)}\,f(k(\lambda)).
$$
Conjugation with $\beta$ gives an automorphism $\alpha$ of $\fA$ since
one verifies that
$$
\alpha( U_\lambda):=\beta U_\lambda\beta^{-1}=D(\lambda) U_{k(\lambda)}.
\ref{\alphaconstr}
$$
In fact,
$$\eqalign{
&\beta U_\lambda\beta^{-1}\phi_\mu=
{1\over D\left(k^{-1}(\mu)\right)}\,\beta U_\lambda
\phi_{k^{-1}(\mu)}={r\left(\lambda,k^{-1}(\mu)\right)
\over D\left(k^{-1}(\mu)\right)}\,\beta\phi_{\lambda*k^{-1}(\mu)}\cr 
&={r\left(\lambda,k^{-1}(\mu)\right)D(\lambda*k^{-1}(\mu))
\over D\left(k^{-1}(\mu)\right)}\,\phi_{k\left(\lambda*k^{-1}(\mu)\right)}
={r\left(\lambda,k^{-1}(\mu)\right)D(\lambda*k^{-1}(\mu))
\over D\left(k^{-1}(\mu)\right)}\,\phi_{k(\lambda)*\mu}.
\cr}$$
Notice that by \udefref\ we have
$$
U_{k(\lambda)}\phi_\mu=r(k(\lambda),\mu)\,\phi_{k(\lambda)*\mu}.
$$
Consequently,
$$
\beta U_\lambda\beta^{-1}\phi_\mu=
{r\left(\lambda,k^{-1}(\mu)\right)D(\lambda*k^{-1}(\mu))
\over D\left(k^{-1}(\mu)\right)r(k(\lambda),\mu)}\,
U_{k(\lambda)}\phi_\mu=D(\lambda) U_{k(\lambda)}\phi_\mu
$$
by \rnormref.

Define 
$$\tau(A):=(\phi_1,A\phi_1)$$
Since $\beta\phi_1=\phi_1$, the state $\tau$ is $\alpha$ invariant.
Moreover vector $\phi_1$ is cyclic and separating for $\fA$ and so the GNS Hilbert
space of state $\tau$ is canonically identified with $\cK$. In this
identification $U_\lambda$ is mapped to $\phi_\lambda$ and the unitary
operator in $\cK$ defined by $\alpha$ is simply $\beta$. Also
$\tau$ is a trace since $\fA$ is linearly generated by $U_\lambda$'s.

We need to verify that the system $\left( \fA,\bZ, \alpha ,\tau \right) $ is
totally ergodic and that its quantum quasi-spectrum 
$(H_1(\alpha),H_2(\alpha),[r_\alpha],k_\alpha)$
is isomorphic to
$(H_1,H_2,[r],k)$. It follows from \betaform\ that the spectrum of
$\beta$ (and equivalently of $\alpha$)
is $H_1$ with $\phi_\lambda$, $\lambda\in H_1$ being the corresponding
eigenvectors. Also
$$
\beta^n f(\lambda)=D\left(\lambda\right)^n
f\left(R(\lambda)^{-n}*\lambda\right),$$
where $R(\lambda):=k(\lambda)*\lambda^{-1}\in H_1$. As $H_1$ has no nontrivial
elements of finite order, $\phi_1$ is the only invariant vector for
$\beta^n$ and $\alpha$ is totally ergodic. Next observe that
$$
\alpha(U_\lambda)={D(\lambda)U_{R(\lambda)}\over r(R(\lambda),\lambda)}
\,U_\lambda,
$$
and so that $H_2(\alpha)$ consists of the operators of the form
${D(\lambda)U_{R(\lambda)}\over r(R(\lambda),\lambda)}$. They are different
for different $\lambda$'s as they correspond to different quasi-eigenvectors
of an ergodic system. The map
$$
H_2\ni\lambda\mapsto\phi(\lambda):={D(\lambda)U_{R(\lambda)}\over
r(R(\lambda),\lambda)}\in H_2(\alpha)
$$
is consequently bijective. $\phi$ is a homomorphism as a consequence of the
following calculation:
$$\eqalign{
\phi(\lambda*\mu)U_{\lambda*\mu}&=\alpha(U_{\lambda*\mu})=
{1\over r(\lambda,\mu)}\alpha(U_\lambda U_\mu)={1\over r(\lambda,\mu)}
\alpha(U_\lambda)\alpha( U_\mu)\cr
&={1\over r(\lambda,\mu)}\phi(\lambda)U_\lambda\phi(\mu)U_\mu=
{1\over r(\lambda,\mu)}\phi(\lambda)*\phi(\mu)U_\lambda U_\mu=
\phi(\lambda)*\phi(\mu)U_{\lambda*\mu}
\cr}$$
Since $R(\phi(\lambda))=R(\lambda)$, it follows that $k=\phi^{-1}k_\alpha\phi$.
Finally, as $\{\phi_\lambda\}$ is a normalized basis in $\cK$ consisting of 
quasi-eigenvectors of $\alpha$, formula \udefref\ implies that
$[r]=\phi^*[r_\alpha]$.
$\square$

\section\examplsec{Examples: Quantum Torus}
\medskip

In this section we consider examples of  systems, defined on  on quantum tori,
illustrating our theory. The first example is a system
satisfying all the assumptions of our classification scheme. Interestingly, it appears as a quantization of a kicked rotor in [BB].

Recall that the algebra 
$\fA$ of observables on a quantum torus is defined as the universal von
Neumann algebra generated by two unitary generators $U,V$ satisfying 
the relation [R]:
$$
UV=e^{2\pi ih}VU\ .
$$
One can think of the elements of $\fA$ as series of the form 
$a=\sum a_{n,m}U^nV^m$. A natural trace on $\fA$ is simply given by
$\tau(a)=a_{0,0}$. The automorphism $\alpha$ is defined on generators by:
$$
\alpha(U):=e^{2\pi i\omega}U,\ \ 
\alpha(V):=UV.
$$
It extends to an automorphism of $\fA$. If $\omega$ is irrational,
then $\alpha$ is totally ergodic. In fact, the eigenvectors of $\alpha$
are just powers of $U$:
$$
\alpha(U^n)=e^{2\pi in\omega}U^n.
$$
Consequently $H_1=\{e^{2\pi in\omega},\ n\in\bZ\}\cong\bZ$ and the spectrum is
simple which proves total ergodicity if $\omega$ is irrational.
Moreover
$$
\alpha(U^nV^m)=e^{2\pi i(n\omega+hm(m-1)/2)}U^m\cdot U^nV^m,
$$
which shows that $U^nV^m$ are quasi-eigenvectors of the second order for $\alpha$.
Since they form an orthonormal basis in $L^{2}\left( \fA,\tau \right) $
we see that $\left( \fA,\bZ, \alpha ,\tau \right) $ is a totally
ergodic system with purely quasi-discrete spectrum of the second order.

We can identify $H_2\cong \bZ^2$ as groups and $H_1$ is simply
the subgroup $\bZ\times\{0\}\subset \bZ^2$. The mapping $R$ is given by
$$
R(n,m)=m\in\bZ\cong H_1,
$$
and the isomorphism $k$ is 
$$k(n,m)=(n+m,m).$$
Define
$$
C(n,m):=\cases{1&if $m=0$\cr
e^{\pi i\left(h(nm-n)+\omega(n^2/m-n)\right)}&otherwise.\cr
}
$$
Then a simple calculation shows that $C(n,m)U^nV^m$ is a normalized basis
for this ergodic system.

In this simple example the group $H_2$ is abelian. 
We can identify $H^2_k(H_2)$ with $H^2(H_2)$,
the second cohomology group of $H_2$. The later group is identified
with the set of symplectic bicharacters by \rsigmarel. A simple calculation shows
that the following symplectic bicharacter represents $[r]$ in our example.
$$
\sigma\left((n,m),(n',m')\right)=e^{2\pi ih(nm'-n'm)},
$$
Notice that $\sigma$ is trivial on $H_1$ and k-invariant.

The above example can be easily extended to give systems with a basis
consisting of quasi-eigenvectors of arbitrary order. Here is one way to
do it. Consider the algebra as before but with an extra generator
$W$ which we assume for simplicity to commute with $U$ and $V$. As before
define a trace $\tau$ such that $\tau(U^nV^mW^k)=0$ unless
$n=m=k=0$. Finally extend the automorphism $\alpha$ by
$\alpha(W)=UVW$. Then one easily verifies that quasi-eigenvectors of
order one are powers of $U$, quasi-eigenvectors of the second order are
$U^nV^m$ and quasi-eigenvectors of the third order are
$U^nV^mW^k$. The last expressions form a basis in the corresponding Hilbert space.

Systems that are ergodic but not totally ergodic are usually associated
with elements of the finite order. For example in an algebra
generated by two unitary generators $U,V$ satisfying 
the relations $UV=e^{2\pi i/N}VU$, and $V^N=1$ consider an automorphism
$\alpha$ given by
$$
\alpha(U):=e^{2\pi i\omega}U,\ \ 
\alpha(V):=e^{2\pi i/N}V.
$$
Here $N$ is a positive integer and $\omega$ is assumed to be irrational.
The eigenvectors of $\alpha$ are $U^nV^m$, $0\leq m\leq N-1$. This
system is ergodic but not totally ergodic since $\alpha^N(V^m)=V^m$
for any $m$.

\bigskip\noindent
{\bf Acknowledgment.} I would like to thank Andrzej Lesniewski for
inspiration, Jerry Kaminker for his enthusiasm and an anonymous referee
for invaluable suggestions.

\vfill\eject

\centerline{\bf References}
\baselineskip=12pt
\frenchspacing

\bigskip

\item{[Ab]} Abramov, L. M.: Metric automorphisms with quasi-discrete spectrum,
{\it Izv. Akad. Nauk USSR}, {\bf 26}, 513--530 (1962)

\item{[Ar]} Araki, H.: Positive cones for von Neumann algebras, 
in Operator Algebras and its Applications, {\it Proc. of Symp. in Pure Math of
the AMS}, vol 38, AMS (1982)

\item{[BB]} J. Bellissard, J., Barelli, A.:
{\it Dynamical localization: a Mathematical Framework}
in Quantum Chaos, Quantum Measurements, edited by P. Cvitanovic, I.C. Percival \& A. Wirzba, Kluwer Publ. (1992)

\item{[B]} Brown, K. S.: {\it Cohomology of Groups}, Springer Verlag (1982)

\item{[BR]}  Bratteli, O., and Robinson, D.: {\it Operator Algebras
and Quantum Statistical Mechanics}, Vol. 1, Springer Verlag (1979)

\item{[C]}  Connes, A.: {\it Noncommutative Geometry}, Academic Press
(1994)

\item{[CFS]}  Cornfeld, I. P., Fomin, S. V., Sinai, Ya.: {\it Ergodic
Theory}, Springer Verlag (1982)

\item{[G]}  Greenleaf, F. P.: {\it Invariant Means on Topological
Groups and Their Applications,} Van Nostrand Reinhold Co. (1969)

\item{[HaP]} Hahn, F., Parry, W.: Minimal dynamical systems with quasi-discrete
spectrum, {\it J. London Math. Soc.}, {\bf 40}, 309--323 (1965)

\item{[H]} Halmos, P.: Measurable transformations, {\it Bull. AMS}, {\bf 55},
1015--1034 (1949)

\item{[HvN]}  Halmos, P., von Neumann, J.: Operator methods in
classical mechanics, II, {\it Ann. Math.}, {\bf \ 43}, 332 - 350 (1942)

\item{[HoP]} Hoare, H., Parry, W.: Affine transformations with quasi-discrete
spectrum I, II, {\it J. London Math. Soc.}, {\bf 41}, 88--96 and 529--530 (1966)

\item{[HR]}  Hewitt, E., Ross, K. A.: {\it Abstract Harmonic Analysis}
, Vol. I, Springer Verlag (1963)

\item{[J]}  Jajte, R.: {\it Strong Limit Theorems in Non-Commutative
Probability, Lecture Notes in Mathematics, }{\bf 1110}, Springer Verlag
(1985)

\item{[K]} Koopman, B.: Hamiltonian Systems and Transformations in Hilbert Spaces,
{\it Proc. Nat. Acad. Sci.}, {\bf 17}, 315--318 (1931)


\item{[KL1]} Klimek, S., Le\'{s}niewski, A.: Quantum Maps, in {\it Contemporary
Mathematics}, {\bf 214}, Proceedings of the 1996 Joint Summer Research Conference
on Quantization, edited by L. Coburn and M. Rieffel, AMS (1997)

\item{[KL2]} Klimek, S., Le\'{s}niewski, A.: Quantum Ergodic Theorems, in 
{\it Contemporary Mathematics}, {\bf 214}, Proceedings of the 1996 Joint Summer
Research Conference on Quantization, edited by L. Coburn and M. Rieffel, AMS
(1997)

\item{[KLMR]}  Klimek, S., Le\'{s}niewski, A., Maitra, N., Rubin, R.:
Ergodic properties of quantized toral automorphisms, 
{\it J. Math. Phys},  {\bf 38}, 67--83 (1997)

\item{[L]}  Lance, C. E.: Ergodic theorems for convex sets and operator
algebras, {\it Inv. Math.,}{\bf \ 37}, 201 - 214 (1976)

\item{[OPT]} Olsen, D., Pedersen, G., Takesaki, M.:
Ergodic actions of compact abelian groups, {\it J. Oper. Theory}, {\bf 3},
237--269 (1980)

\item{[R]}  Rieffel, M.: Non-commutative tori - a case study of
non-commutative differentiable manifolds, {\it Contemp. Math}., {\bf 105}
, 191 -- 211 (1990)

\item{[Si]}  Sinai, Ya.: {\it Topics in Ergodic Theory}, Princeton
University Press (1994)

\item{[VN]}  Von Neumann, J.: Zur Operathorenmethode in der klassischen
Mechanik, {\it Ann. Math., }{\bf 33}, 587 - 642 (1932)

\item{[W]}  Walters, P.: {\it An Introduction to Ergodic Theory},
Springer Verlag (1982)

\end